\documentclass[12pt]{article}
= 0.0in

\def\Dj{\hbox{D\kern-.73em\raise.30ex\hbox{-}
\raise-.30ex\hbox{}}}
\def\dj{\hbox{d\kern-.33em\raise.80ex\hbox{-}
\raise-.80ex\hbox{\kern-.40em}}}
\usepackage{amsmath,amsthm,amsfonts,amssymb,amscd,cite,graphicx}
\allowdisplaybreaks

\newtheorem{theorem}{Theorem}
\newtheorem{corollary}{Corollary}
\newtheorem{conjecture}{Conjecture}
\newtheorem{lemma}{Lemma}
\newtheorem{definition}{Definition}

\newtheorem{note}{Note}

\begin{document}

\baselineskip=0.30in

\vspace*{2cm}

\noindent {\bf \Large On Minimum Terminal Distance Spectral Radius\\
of Trees with Given Degree Sequence}

\vspace{6mm}

\noindent {\bf \large Mikhail Goubko}\footnote{Corresponding author
at: V.A. Trapeznikov Institute of Control
Sciences,\hspace*{5mm}Russian Academy of Science,
\hspace*{5mm} Moscow, Russia. \\
\hspace*{5mm} {\it E-mail addresses:} {\tt mgoubko@mail.ru} (M.
Goubko)}

\vspace{5mm}

\baselineskip=0.20in

\noindent {\it V.A. Trapeznikov Institute of Control Sciences,\\
Russian Academy of Science, Moscow, Russia}

\vspace{5mm}

\noindent
{\it Keywords}\\
Terminal distance spectral radius\\
Extremal tree with given degree sequence\\
Wiener index for vertex-weighted graphs \\

\vspace{4mm}


\vspace{10mm}

\noindent {\bf A \ B \ S \ T \ R \ A \ C \ T}

\vspace{2mm}

\noindent

For a tree with the given sequence of vertex degrees the spectral
radius of its terminal distance matrix is shown to be bounded from
below by the the average row sum of the terminal distance matrix of
the, so called, $BFS$-tree (also known as a \emph{greedy tree}).
This lower bound is typically not tight; nevertheless, since
spectral radius of the terminal distance matrix of $BFS$-tree is a
natural upper bound, the numeric simulation shows that relative gap
between the upper and the lower bound does not exceed $3\%$ (we also
 make a step towards justifying this fact analytically.) Therefore,
 the conjecture that $BFS$-tree has the minimum terminal distance
 spectral radius among all trees with the given degree sequence is valid
 with accuracy at least $97\%$. The same technique can be applied to the
 \emph{distance spectral radius} of trees, which is a more popular
 topological index.
 \vspace{10mm}

\baselineskip=0.30in

\section{Introduction}

We study simple connected graphs. Let $V(G)$ be the vertex set and
$E(G)$ be the edge set of an undirected graph $G$. For any pair of
vertices $u,v\in V(G)$ let $d_G(u, v)$ denote the distance (the
length of the shortest path) between $u$ and $v$ in $G$.

The matrix $D(G):=(d_G(u,v))_{u,v\in V(G)}$ is known as a
\emph{distance matrix} of a graph $G$. Along with the
\emph{adjacency matrix} and the {\it Laplacian matrix\/}
\cite{mohar}, the $D(G)$ and related matrices are the most popular
objects to study in algebraic graph theory.

Denote by $d_G(v)$ the degree of a vertex $v \in V(G)$ in the graph
$G$, i.e., the number of vertices being adjacent to $v$. The vertex
$v \in V(G)$ is said to be {\it pendent\/} if $d_G(v)=1$. All other
vertices of the graph $G$ are referred to as {\it internal\/}. By
$W(G)$ we denote the set of all pendent vertices of a graph $G$, and
let $M(G) = V(G) \setminus W(G)$ be the set of its internal
vertices.

The matrix $RD(G)=(d_G(u,v))_{u,v\in W(G)}$ is typically referred to
as the \emph{terminal distance matrix} or the \emph{reduced distance
matrix} of graph $G$. A \emph{tree} is a connected graph with $N$
vertices and $N-1$ edges. Terminal distance matrices of trees are of
special interest, since a tree can be reconstructed by its terminal
distance matrix (see \cite{Zaretskii1965}).

Concepts based on the distance matrix are intensively employed in
the mathematical chemistry. In particular, one of the oldest
topological molecular indices, the \emph{Wiener index}, is defined
as one half of the sum of all elements of the distance matrix of a
graph:
$$WI(G)=\frac{1}{2}\sum_{u,v\in V(G)}d_G(u,v),$$
and represents just an example from the large family of
\emph{distance-based topological indices} \cite{Todeschini2000,
Todeschini2009}. The \emph{terminal Wiener index} of graph $G$ is
defined by analogy as the one half of a sum of elements of $RD(G)$:
$$TWI(G)=\frac{1}{2}\sum_{u,v\in W(G)}d_G(u,v).$$

\emph{Spectrum-based indices}, which are calculated using
eigenvalues and eigenvectors of various graph matrices, form a yet
another family of topological indices \cite{Balaban1991}, the most
famous being the \emph{Estrada index} \cite{Estrada2000}. Balaban et
al. \cite{Balaban1991} suggested the \emph{distance spectral radius}
(the largest eigenvalue of the distance matrix) as a molecular
descriptor giving rise to the extensive QSPR\footnote{QSPR =
quantitative structure-property relations.} research and to the
studies of mathematical properties of the distance spectral radius
($DSR$).

In particular, in \cite{Gutman1998} upper and lower bounds were
suggested (and later improved by Zhou \cite{Zhou2007}) for $DSR$ of
a tree in terms of the tree order and the sum of squares of
distances between unordered pairs of vertices. Das \cite{Das2009}
obtained bounds for $DSR$ for bipartite graphs in terms of partition
orders and characterized extremal graphs.

Indulal \cite{Indulal2009} has shown that
\begin{equation}\label{eq_DSRgeWI}
DSR(G)\ge \frac{2}{n}WI(G)
\end{equation}
for any
connected graph $G$ with $n$ vertices, with equality if and only if
$G$ is \emph{distance regular}, i.e., when all row sums
$D_1,...,D_n$ of $D(G)$ (also known as \emph{distance degrees}) are
equal to each other. He also proved the inequality $DSR(G)\ge
\sqrt{\frac{D_1^2+...+D_n^2}{n}}$ with equality if and only if $G$
is distance regular, and even stronger lower bounds in terms of row
sums, the \emph{second degree distance sequence}, and the Wiener
index (see \cite{Indulal2009} for details). Alternative bounds of
this sort were reported later by He et al \cite{He2010}. In
\cite{ZhouIlic2010} a simple lower bound has been suggested in terms
of graph order and of two maximal vertex degrees. Upper and lower
bounds of $DSR$ were also obtained in \cite{ZhouIlic2010} for
bipartite graphs in terms of graph order, diameter, and extremal
degrees in two partitions. Recent results are summarized in
\cite{StevanovicIlic2012}.

\begin{definition} \emph{\cite{Marshall79,Zhang2008W}}
For a real sequence $\mathbf{x}=(x_1,..., x_n)$, $n\in \mathbb{N}$,
denote with $\mathbf{x}^\uparrow=(x_{[1]},..., x_{[n]})$ the
sequence, where all components of $\mathbf{x}$ are arranged in
ascending order.
\end{definition}

For a graph $G$ define its \emph{degree sequence} as
$\mathbf{d}(G):=(d_G(v))_{v\in V(G)}^\uparrow$.

\begin{definition}A sequence $\mathbf{d}=(d_1, ..., d_n)$ is called \emph{graphic}
if such a graph $G$ exists that $\mathbf{d}(G)=\mathbf{d}$
\emph{(}see \emph{\cite{Zhang2008W}}\emph{)}. If $G$ is a tree, we
will say that $\mathbf{d}$ \emph{generates a tree}.
\end{definition}

A non-decreasing natural sequence $\mathbf{d}=(d_1, ..., d_n)$ is
known to generate a tree if and only if $d_1+...+d_n=2(n-1)$. Let
$\mathcal{T}(\mathbf{d})$ be the set of all trees with the degree
sequence $\mathbf{d}$.

An upper bound for $DSR(G)$ has been suggested in
\cite{ZhouIlic2010} in terms of minimum degree, second minimum
degree and the diameter of graph $G$. Later in \cite{Lin2013,
Chen2013} a series of improved upper bounds were formulated for
$DSR(G)$ in terms of degree sequence and distance matrix row sums of
$G$.

In \cite{Zhang2008L} Zhang introduced the notion of
\emph{$BFS$-tree} (breadth-first-search tree) and proved that
$BFS$-tree has the largest spectral radius of the Laplacian matrix
among all trees with the given degree sequence. The $BFS$-tree
generated by a degree sequence $\mathbf{d}=(d_1,...,d_n)$ is denoted
as $BFS(\mathbf{d})$ and is built in a ``top-down'' manner by
sequentially filling ``levels''. We start with a single vertex at
Level 1 and connect it to $d_n$ vertices at Level 2. Then we add
$d_{n-1}+...+d_{n-d_n}-d_n$ vertices to Level 3 and connect them
sequentially ``left to right'' to $d_n$ vertices from Level 2 so
that the latter have degrees $d_{n-1},...,d_{n-d_n}$ respectively.
In the same manner we add vertices to Level 4 and connect them
sequentially ``left to right'' to the vertices from Level 3 picking
the largest unused degree from the sequence $\mathbf{d}$, and so
forth until all $n$ vertices are added to the tree (see the example
in Fig.~\ref{figureBFS} and \cite{Zhang2008L, Zhang2008W} for the
detailed algorithm). Later $BFS$-tree (also known as a \emph{greedy
tree}) was shown to minimize the Wiener index \cite{Zhang2008W,
Wang2008}, the terminal Wiener index \cite{Szekely2011}, and the
maximal number of subtrees \cite{ZhangWangSubtrees2013} among all
trees with the given degree sequence.

\vspace{4mm}

\begin{center}
\includegraphics[height=4cm,keepaspectratio]{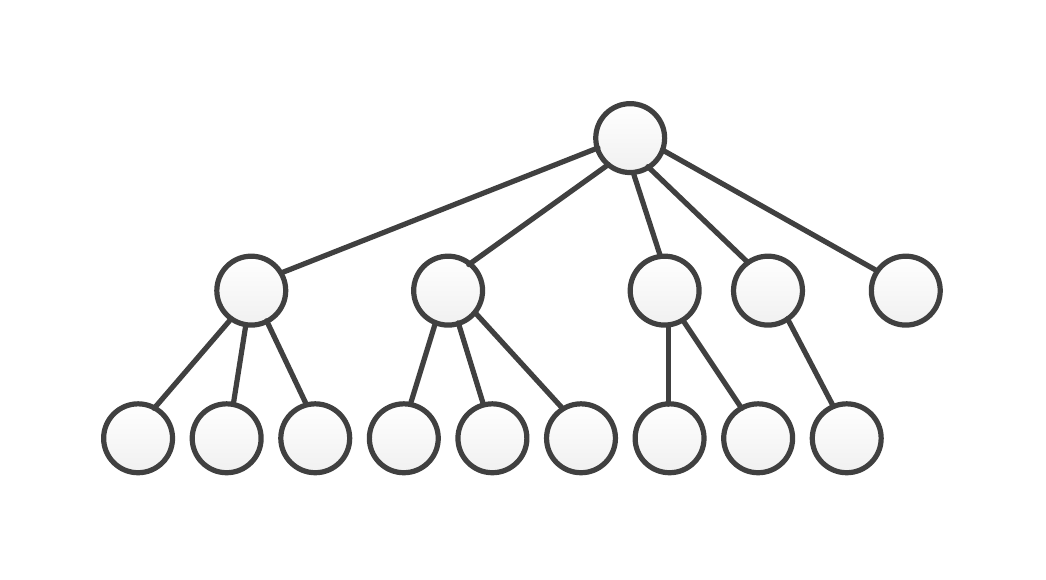}
\label{figureBFS}
\end{center}
\baselineskip=0.20in

\noindent {\bf Fig. 1.} An example of $BFS$-tree for the degree
sequence $(\underbrace{1,...,1}_{10\textrm{ times}},2,3,4,4,5)$.

\baselineskip=0.30in

\vspace{3mm}

A special case of the $BFS$-tree for the degree sequence
$(1,...,1,r, \Delta,..., \Delta)$, where $1 \le r \le \Delta$, is
referred to as the \emph{Volkmann tree}. It is shown
\cite{WagnerZhangWang2013} to be extremal with respect to the Wiener
index and to many other distance-based topological indices over all
trees of order $n$ and maximum degree $\Delta$.

Stevanovic and Ilic \cite{StevanovicIlic2010} conjectured that the
Volkmann tree has the minimum $DSR$ among all trees of order $n$ and
maximum degree $\Delta$. Taking into account the above
considerations, this conjecture can be generalized as follows:

\begin{conjecture}\label{conj1}
The $BFS$-tree has the minimum distance spectral radius among all
trees with the given degree sequence.
\end{conjecture}

The \emph{terminal distance spectral radius} ($TDSR$) is defined by
analogy to $DSR$ as the spectral radius of the terminal distance
matrix $RD(G)$ of graph $G$. This index is less studied in the
literature, yet its behavior is alike that of $DSR$ (at least, for
trees). Therefore, a conjecture similar to Conjecture \ref{conj1}
can be put for $TDSR$:

\begin{conjecture}\label{conj2}
The $BFS$-tree has the minimum terminal distance spectral radius
among all trees with the given degree sequence.
\end{conjecture}

In the present paper we do not prove these conjectures in full.
Nevertheless, below we show that $DSR$ of the $BFS$-tree is at least
very close to the minimum $DSR$ of a tree with the given degree
sequence, and $TDSR$ of the $BFS$-tree is very close to the minimum
$TDSR$.

More formally, we show in Section \ref{section_simple} that $DSR$ of
a tree with degree sequence $\mathbf{d}$ having $n$ pendent vertices
is bounded from below by $TLB(\mathbf{d})$, the average row sum of
the distance matrix of the $BFS$-tree (or, equivalently, by the
value $\frac{2}{n} WI(BFS(\mathbf{d}))$). On the other hand, the
minimal value of $DSR$ does not exceeds $TUB(\mathbf{d})$, which is
defined as the distance spectral radius of the $BFS$-tree. That is,
if a tree $T^*$ minimizes $DSR$ over $\mathcal{T}(\mathbf{d})$, we
show that
\begin{equation}\label{eq_main_ineq_dsr}
LB(\mathbf{d}):=\frac{2}{n} WI(BFS(\mathbf{d}))\le DSR(T^*)\le
UB(\mathbf{d}):=DSR(BFS(\mathbf{d})).
\end{equation}
Similar inequalities trivially hold for the tree $T^{**}$, which has
the minimum terminal distance spectral radius:
\begin{equation}\label{eq_main_ineq_tdsr}
TLB(\mathbf{d}):=\frac{2}{n} TWI(BFS(\mathbf{d}))\le TDSR(T^{**})\le
TUB(\mathbf{d}):=TDSR(BFS(\mathbf{d})).
\end{equation}

After some preliminary definitions in Section
\ref{section_preliminaries} we provide the alternative justification
of inequalities (\ref{eq_main_ineq_dsr}) and
(\ref{eq_main_ineq_tdsr}) in terms of \emph{the Wiener index for
vertex-weighted graphs}, the distance-based index first introduced
in \cite{Klavzar1997} and later studied in \cite{Goubko2015,
GoubkoMilos2015}.

In Section \ref{section_quality} we show that $TLB(\mathbf{d})$ is
very close to $TUB(\mathbf{d})$, and, therefore,
$TDSR(BFS(\mathbf{d}))$ approximates the extremal value
$TDSR(T^{**})$ well.

In the conclusion we outline open issues and perspectives.

\section{Main inequalities} \label{section_simple}

For a symmetric real $n\times n$ matrix $A$ denote with
$\lambda_1(A), ..., \lambda_n(A)$ its eigenvalues enumerated in
descending order. The biggest eigenvalue of matrix $A$ is called its
\emph{spectral radius}. The corresponding eigenvector is called
\emph{Perron vector} and is denoted as $\mathbf{p}(A)$. Let $S_n$ be
the unit sphere in $\mathbb{R}_n$, and let $S_n^+$ be its
intersection with the non-negative orthant.

Since, by inequality (\ref{eq_DSRgeWI}), $DSR(G)\ge
\frac{2}{n}WI(G)$ for any connected graph of order $n$, we conclude
immediately that for any degree sequence $\mathbf{d}$ generating a
tree with $n$ vertices and any $T\in\mathcal{T}(\mathbf{d})$ we have
$$DSR(T)\ge
\min_{G\in\mathcal{T}(\mathbf{d})} \frac{2}{n}WI(G).$$

From \cite{Wang2008, Zhang2008W} we know that the Wiener index is
minimized over $\mathcal{T}(\mathbf{d})$ by the corresponding
$BFS$-tree. On the other hand, since
$BFS(\mathbf{d})\in\mathcal{T}(\mathbf{d})$, the minimum value of
$DSR$ trivially does not exceed $DSR(BFS(\mathbf{d}))$, and
inequalities (\ref{eq_main_ineq_dsr}) follow immediately.

For $TDSR$ an inequality similar to (\ref{eq_DSRgeWI}) is easily
obtained. Since, by the Raileigh-Ritz principle,
$\lambda_1(A)=\max_{\mu\in S_n}\mathbf{\mu}^T A \mathbf{\mu}$ for
any real symmetric $n\times n$ matrix $A$, we conclude that for any
connected graph $G$ with $n$ pendent vertices
$$\lambda_1(RD(G))\ge \left(\frac{1}{\sqrt{n}},...,\frac{1}{\sqrt{n}}\right) RD(G) \left(\frac{1}{\sqrt{n}},...,\frac{1}{\sqrt{n}}\right)^T\ge$$
$$\ge \frac{1}{n}\left(1,...,1\right) RD(G) \left(1,...,1\right)^T=\frac{2}{n}TWI(G),$$
and inequalities (\ref{eq_main_ineq_tdsr}) are obtained similar to
(\ref{eq_main_ineq_dsr}).

\section{Huffman Trees and Generalized Wiener Index} \label{section_preliminaries}

Below we provide an alternative proof of inequalities
(\ref{eq_main_ineq_dsr}) and (\ref{eq_main_ineq_tdsr}) using recent
results \cite{Goubko2015,GoubkoMilos2015} on minimization of the
Wiener index for vertex-weighted trees. In this section we recall
basic definitions and theorems following \cite{GoubkoMilos2015}.

A graph $G$ is called \emph{vertex-weighted}, if each vertex $v\in
V(G)$ is endowed with a non-negative weight $\mu_G(v)$.
%With $\mu_G$ we denote the total vertex weight of the graph $G$, and
%$\mathcal{WT}(n)$ stands for the set of all vertex-weighted trees of
%order $n$.

The \emph{Wiener index for vertex-weighted graphs} is defined in
\cite{Klavzar1997} as
$$VWWI(G):=\frac{1}{2}\sum_{u,v\in V(G)}\mu_G(u)\mu_G(v) d_G(u,v).$$

\begin{definition}
Consider a vertex set $V$ with $|V|=n$. Let the function $\mu: V
\rightarrow \mathbb{R}_+$ assign a non-negative weight $\mu(v)$ to
each vertex $v\in V$, while the function $d: V \rightarrow
\mathbb{N}$ assign a natural degree $d(v)$. The tuple $\langle\mu,
d\rangle$ is called \emph{a generating tuple} if $\sum_{v\in V}
d(v)=2(n-1)$. Denote with $\overline{\mu}:=\sum_{v\in V}\mu(v)$ the
total weight of vertex set $V$.
\end{definition}

Below we use a degree sequence $\mathbf{d}$ as a synonym of a degree
function $d(\cdot)$ using natural correspondence between these
concepts.

Let $\mathcal{WT}(\mu, d)$ be the set of trees over the vertex set
$V$ with vertex weights $\mu(\cdot)$ and degrees $d(\cdot)$, and let
$V(\mu, d)$ be the domain of functions of a generating tuple
$\langle\mu, d\rangle$. Introduce the set $W(\mu, d):=\{w \in V(\mu,
d): d(w)=1\}$ of \emph{pendent} vertices and the set $M(\mu,
d):=V(\mu,d)\backslash W(\mu,d)$ of \emph{internal} vertices.

\begin{definition}
We will say that in a generating tuple $\langle\mu, d\rangle$
\emph{weights are degree-monotone}, if for any $m, m' \in M(\mu,d)$
from $d(m) < d(m')$ it follows that $\mu(m) \le \mu(m')$, and, also,
$\mu(w)>0$ for any $w\in W(\mu,d)$.
\end{definition}

For a generating tuple $\langle\mu,d\rangle$ the \emph{generalized
Huffman algorithm} \cite{Goubko2015} builds a tree $H \in
\mathcal{WT}(\mu,d)$ as follows.

\textbf{Setup.} Define the vertex set $V_1 := V(\mu,d)$ and the
functions $\mu^1$ and $d^1$, which endow its vertices with weights
$\mu^1(v) := \mu(v)$ and degrees $d^1(v) := d(v)$, $v \in V_1$. We
start with the empty graph $H$ over the vertex set $V(\mu,d)$.

\textbf{Steps $i = 1, ..., q-1$.} Denote with $m_i$ the vertex
having the least degree among the vertices of the least weight in
$M(\mu^i,d^i)$. Let $w_1, ..., w_{d(m_i)-1}$ be the vertices having
$d(m_i)-1$ least weights in $W(\mu^i,d^i)$. Add to $H$ edges
$w_1m_i, ..., w_{d(m_i)-1}m_i$.

Define the set $V_{i+1} := V_i \backslash \{w_1, ...,
w_{d(m_i)-1}\}$ and functions $\mu^{i+1}(\cdot), d^{i+1}(\cdot)$,
endowing its elements with weights and degrees as follows:
\begin{gather}
\mu^{i+1}(v) := \mu^i(v)\text{ for }v \neq m_i, \hspace{20pt} \mu^{i+1}(m_i) := \mu^i(m_i)+\mu^i(w_1)+...+\mu^i(w_{d(m_i)-1}),\nonumber\\
\label{eq_Huffman_tuples}d^{i+1}(v) := d^i(v)\text{ for }v \neq m_i,
\hspace{20pt} d^{i+1}(m_i) := 1.
\end{gather}

\textbf{Step $q$.} Consider a vertex $m_q \in M(\mu^q, d^q)$. By
construction, $|M(\mu^q, d^q)|=1$, $|W(\mu^q, d^q)|=d(m_q)$. Add to
$H$ edges connecting all vertices from $W(\mu^q, d^q)$ to $m_q$.
Finally, set $\mu_H(v) := \mu(v)$, $v \in V(H)$.

In general, the Huffman tree is not unique, since there can be more
than one set of vertices having $d(m_i)-1$ least weights in
$W(\mu^i,d^i)$ at Step $i$. Denote with $\mathcal{H}(\mu, d)$ the
set of all Huffman trees for the generating tuple
$\langle\mu,d\rangle$.

\begin{theorem}\label{theorem_huffman_vwwi}\textbf{\emph{\cite{Goubko2015}}} \sloppy If weights are degree-monotone
in a generating tuple $\langle\mu,d\rangle$, then a tree $T$
minimizes the Wiener index over the set $\mathcal{WT}(\mu,d)$ of
trees with given vertex weights and degrees, if and only if
$T\in\mathcal{H}(\mu, d)$.
\end{theorem}

If $\mathbf{\mu}(G):=(\mathbf{\mu}_G(v))_{v\in V(G)}$ is a vector of
vertex weights in graph $G$, $VWWI(G)$ reduces to the quadratic form
$\frac{1}{2}\mathbf{\mu}(G)^TD(G)\mathbf{\mu}(G)$.\footnote{To
simplify notation we always assume below that vector components and
matrix rows go in the same order and no confusion arises.}

The ``classic'' Wiener index $WI$ is a special case of $VWWI$ for
unit weights. It is shown in \cite{Wang2008,Zhang2008W} that for
$T\in \mathcal{T}(\mathbf{d})$ we have $WI(T)\ge
WI(BFS(\mathbf{d}))$ with equality if and only if $T\sim
BFS(\mathbf{d})$. Therefore, from Theorem \ref{theorem_huffman_vwwi}
we conclude that
\begin{note}\label{note_BFSisHuffman1}$BFS(\mathbf{d})$ is isomorphic to some Huffman tree for
unit vertex weights and the same degree sequence $\mathbf{d}$.
\end{note}

If we define the vector of terminal vertices' weights as
$\mathbf{w}(G):=(\mu_G(v))_{v\in W(G)}$, the terminal Wiener index
for vertex-weighted trees is defined by analogy to the ``classic''
terminal Wiener index as
\begin{equation}\label{eq_tvwwi_def}
TVWWI(G)=\frac{1}{2}\sum_{u,v\in
W(G)}\mu_G(u)\mu_G(v)d_G(u,v)=\frac{1}{2}\mathbf{w}(G)^TRD(G)\mathbf{w}(G).
\end{equation}

This index is a special case of $VWWI$ for internal vertices having
zero weights, and the terminal Wiener index $TWI$ is a special case
of $TVWWI$ for terminal vertices having unit weights. Similar to the
Wiener index, it is shown in \cite{Szekely2011} that $TWI(T)\ge
TWI(BFS(\mathbf{d}))$ for all $T\in \mathcal{T}(\mathbf{d})$. Hence,
from Theorem \ref{theorem_huffman_vwwi} we see that
\begin{note}\label{note_BFSisHuffman2}$BFS(\mathbf{d})$ is isomorphic to some Huffman
tree for unit weights of pendent vertices, zero weights of internal
vertices, and degree sequence $\mathbf{d}$.
\end{note}

\begin{note}\label{note_tvwwi_monotone}\sloppy
If in a generating tuple $\langle\mu,d\rangle$ internal vertices
have zero weights, then weights are always degree-monotone in
$\langle\mu,d\rangle$ and, by Theorem
\emph{\ref{theorem_huffman_vwwi}}, $TVWWI(G)$ is minimized over
$\mathcal{WT}(\mu,d)$ with some Huffman tree.
\end{note}

\section{Lower Bounds: Alternative Proofs} \label{section_lb}

First we estimate from below $TDSR(T):=\lambda_1(RD(T))$ of a tree
with the given degree sequence. If a vector
$\mu=(\mu_1,...,\mu_n)^T$ and a tree degree sequence
$\mathbf{d}=(d_1,...,d_N)$ are given, where $n\le N$, denote with
$\langle\mu,d\rangle$ a generating tuple obtained by associating
weights to degrees ''left-to-right'' and assigning zero weights to
the rest of the vertices.

\begin{lemma}\label{lemma_TLB_x}\sloppy
If a degree sequence $\mathbf{d}$ has $n$ elements being equal to
unity, then for any tree $T\in \mathcal{T}(\mathbf{d})$ the
inequality holds
\begin{equation}\label{TDSR_ge_TVWWI}TDSR(T)\ge
2\max_{\mu\in S_n^+}TVWWI(H(\mu, d)),
\end{equation}
where $H(\mu,
d)$ is any Huffman tree from $\mathcal{H}(\mu,d)$.
\begin{proof}
Obviously,
$$TDSR(T)\ge \min_{G\in \mathcal{T}(\mathbf{d})}\lambda_1(RD(G)).$$

By the Raileigh-Ritz principle,
$$\lambda_1(RD(G))=\max_{\mu\in S_n}\mathbf{\mu}^T RD(G) \mathbf{\mu},$$
and the maximum is achieved at $\mathbf{\mu}=\mathbf{p}(RD(G))$.
Since any terminal distance matrix $RD(G)$ is positive, by
Perron-Frobenius theorem, all components of the Perron vector
$\mathbf{p}(RD(G))$ are positive, and we can limit maximization to
$S_n^+$. Therefore,
$$TDSR(T)\ge \min_{G\in \mathcal{T}(\mathbf{d})}\max_{\mu\in S_n^+}\mu^T RD(G) \mu.$$

By the famous \emph{minimax inequality} we only decrease the right
side by changing the order of taking the minimum and the maximum.
Therefore,
$$TDSR(T)\ge \max_{\mu\in S_n^+}\min_{G\in \mathcal{T}(\mathbf{d})}\mu^T RD(G) \mu.$$
From (\ref{eq_tvwwi_def}), we have
$$TDSR(T)\ge 2\max_{\mu\in S_n^+}\min_{G\in
\mathcal{WT}(\mu,d)}TVWWI(G).$$

By Note \ref{note_tvwwi_monotone}, the latter minimum is achieved at
some Huffman tree for the generating tuple $\langle\mu,d\rangle$,
and we obtain (\ref{TDSR_ge_TVWWI}).
%Finally,
%$$TDSR(T)\ge 2 \max_{\mu\in S_n^+}TVWWI(H),\text{ where }H\in \mathcal{H}(\mu,d).$$
\end{proof}
\end{lemma}

From (\ref{TDSR_ge_TVWWI}), the lower bound (
\ref{eq_main_ineq_tdsr}) follows immediately, since, by Note
\ref{note_BFSisHuffman2}, $BFS$-tree is isomorphic to one of Huffman
trees for equal pendent vertex weights, and the weight vector
$\frac{1}{\sqrt{n}}(1, ..., 1)^T$ belongs to $S_n^+$. Nevertheless,
below we prove a somewhat stronger result, which may be useful in
many respects, namely, that the maximum in the right-hand side of
(\ref{TDSR_ge_TVWWI}) is attained when all vertex weights are equal.

\begin{theorem}\label{theorem_bfs_optimal}
If a degree sequence $\mathbf{d}$ generates a tree with $n$ pendent
vertices, then for any vector $\mu\in S_n^+$ and a tree $T\in
\mathcal{WT}(\mu, d)$ the following inequality holds:
\begin{equation}\label{ineq_BFS_worst}
TWI(BFS(\mathbf{d}))\ge n\cdot TVWWI(T).
\end{equation}
\end{theorem}
We postpone the proof to the Appendix. Note that the reasoning there
(namely, the Corollary \ref{cor_simplex}) implies that a similar
inequality $TWI(BFS(\mathbf{d}))\ge n^2 TVWWI(T)$ is valid in an
even more restricting environment when the positive vector $\mu$ of
pendent vertex weights is taken from the $n$-dimensional simplex
$\mu_1+...+\mu_n=1$.

We can repeat the same line of proof for $DSR$. Lemma
\ref{lemma_TLB_x} extends immediately providing an alternative proof
of inequalities (\ref{eq_main_ineq_dsr}).

Extension of Theorem \ref{theorem_bfs_optimal} (i.e., proving that
equal vertex weights in $VWWI$ make the best estimate of $DSR$) can
be useful in general, but in the context of the present paper is a
side line of the analysis and can be skipped.

\section{Quality of Lower Bound}\label{section_quality}

In the previous sections we presented the lower bound
$TLB(\mathbf{d})$ and the upper bound $TUB(\mathbf{d})$ for the
\textbf{minimum} terminal distance spectral radius of a tree with
the given degree sequence $\mathbf{d}$. We have shown that the
terminal distance spectral radius $TDSR(T^{**})$ of such an extremal
tree $T^{**}\in \mathcal{T}(\mathbf{d})$ lies between the average
terminal distance row sum and the terminal distance spectral radius
of the $BFS$-tree (which can be thought as a ``maximally balanced
tree'' for the degree sequence $\mathbf{d}$) .

When the average row sum (or, equivalently, the value of
$\frac{2}{n}TWI(BFS(\mathbf{d}))$) is equal to
$TDSR(BFS(\mathbf{d}))$, Conjecture \ref{conj2} holds, and the
$BFS$-tree minimizes $TDSR$. But, by the Perron-Frobenius theorem,
it is true only when $BFS$-tree is terminal distance regular (i.e.,
all row sums of the terminal distance matrix are equal to each
other), which is not a typical case.

Nevertheless, below we, to some extent, justify that even when the
lower bound is not attained, its relative error
\begin{equation}\label{eq_error}
TErr(\mathbf{d}):=\frac{TUB(\mathbf{d})-TLB(\mathbf{d})}{TUB(\mathbf{d})}
\end{equation}
is small (namely, that it does not exceed $3\%$).

Firstly, we provide some computational evidence. In
Fig.~\ref{figerror} the relative error (\ref{eq_error}) is presented
for all possible degree sequences with the tree order $N$ not
exceeding $22$.

\vspace{4mm}

\begin{center}
\includegraphics[width=6cm,keepaspectratio]{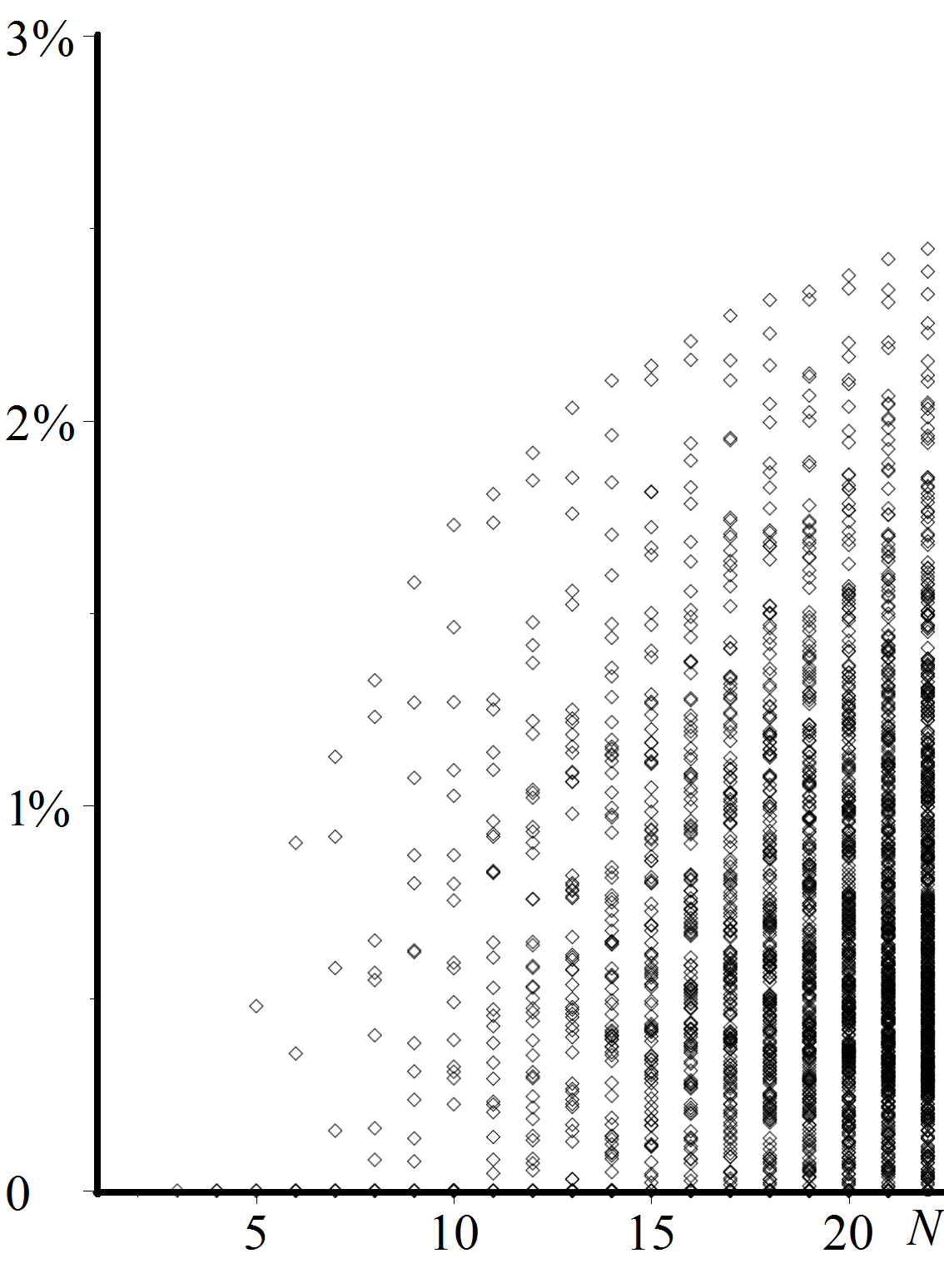}
\label{figerror}
\end{center}
\baselineskip=0.20in

\noindent {\bf Fig. 2.} Relative error of the lower bound for $TDSR$
vs the tree order $N$.

\baselineskip=0.30in

\vspace{3mm}

We see that the relative error never exceeds $3\%$ with most degree
sequences having error less than $1\%$. Analysis of ``extremal''
degree sequences (those constituting the upper envelope of the point
set in the figure), shows that all they have the form
$\mathbf{d}(a,b):=(\underbrace{1,...,1}_{a+b},\underbrace{2,...,2}_{a},a+b)$,
where $a,b\in \mathbb{N}$ (see Fig.~\ref{figBFS34} for the example
of the $BFS$-tree for the degree sequence $\mathbf{d}(3,4)$). This
observation leads us to the following conjecture.

\vspace{4mm}

\begin{center}
\includegraphics[height=3cm,keepaspectratio]{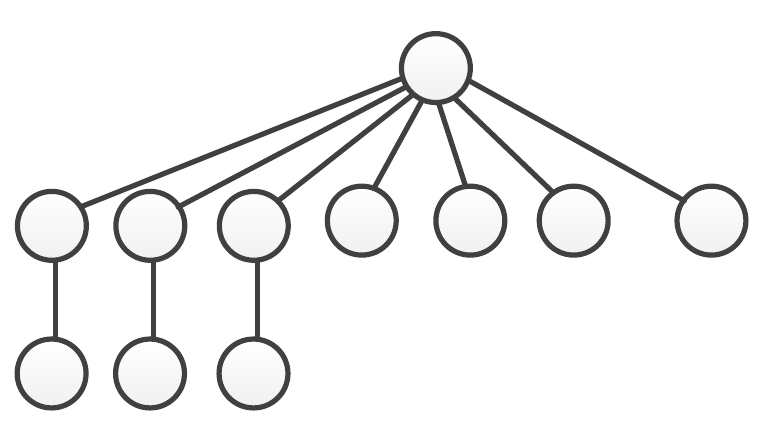}
\label{figBFS34}
\end{center}
\baselineskip=0.20in

\noindent {\bf Fig. 3.} $BFS$-tree for $\mathbf{d}(3,4)$.

\baselineskip=0.30in

\vspace{3mm}

\begin{conjecture}\label{conj3}
For every tree degree sequence $\mathbf{d}$ there exist such $a,b\in
\mathbb{N}$ that $TErr(BFS(\mathbf{d}))\le
TErr(BFS(\mathbf{d}(a,b)))$.
\end{conjecture}

In the following lemma we obtain a closed-form expression for $TDSR$
of the $BFS$-tree for the degree sequence $\mathbf{d}(a,b)$.
\begin{lemma}\label{lemma_TDSR_ab}
\begin{equation}\label{eq_TDSR_ab}
TDSR(BFS(\mathbf{d}(a,b)))=2a+b-3+\sqrt{4a^2+b^2+5ab-4a+2b+1}.
\end{equation}
\begin{proof}
Denote with $I_n$ the identity matrix of order $n$ and let $E_{mn}$
stand for the all-ones $m\times n$ matrix. Consider the terminal
distance matrix $RD$ of the tree $BFS(\mathbf{d}(a,b))$:
\begin{equation}\label{eq_rd}RD=\left(%
\begin{array}{cc}
  4(E_{aa}-I_a)   & 3E_{ab}\\
  3E_{ba}   & 2(E_{bb}-I_b)\\
\end{array}
\right).\end{equation}

Consider a Perron vector
$\mathbf{u}(RD)=\mathbf{x}=(x_1,...,x_{a+b})$. Since any two of the
first $a$ pendent vertices in $BFS(\mathbf{d}(a,b))$ are swapped
with no impact on $RD$, we have $x_1=...=x_a=: \alpha$. Also, the
last $b$ vertices in $BFS(\mathbf{d}(a,b))$ are similar to each
other and, hence, $x_{a+1}=...=x_{a+b}=:\beta$.

As, by definition, $\lambda_1(RD) \mathbf{x}=RD\mathbf{x}$, we
obtain the following system of equations:
$$\left\{%
\begin{array}{l}
    \lambda_1(RD)\alpha=4(a-1)\alpha+3b\beta\\
    \lambda_1(RD)\beta=3a\alpha+2(b-1)\beta\\
\end{array}%
\right.    $$

After elimination of $\alpha$ and $\beta$ we find $\lambda_1(RD)$ as
the positive root of the square equation
$\lambda^2-2(2a+b-3)\lambda-8(a+b-1)-ab=0$, which immediately gives
 (\ref{eq_TDSR_ab}).
\end{proof}
\end{lemma}

To write down $TErr(\mathbf{d}(a,b))$, evaluate
$TLB(\mathbf{d}(a,b))$ as the average row sum of $RD$:
\begin{equation}\label{eq_lb}
TLB(\mathbf{d}(a,b))=\frac{4a(a-1)+2\cdot3ab+2b(b-1)}{a+b}=2(2a+b)\left(1-\frac{1}{a+b}\right).
\end{equation}

Finally, put the following estimate of $TErr(\cdot)$, which is true
if Conjecture \ref{conj3} holds.

\begin{theorem}
If Conjecture \ref{conj3} holds, then for any tree degree sequence
$\mathbf{d}$
$$TErr(\mathbf{d})\le
\frac{3\sqrt{2}-4}{3\sqrt{2}+4}<0.03.$$
\begin{proof}
Using Lemma \ref{lemma_TDSR_ab} and expression (\ref{eq_lb}) write
\begin{equation}\label{eq_TErr_ab}TErr(\mathbf{d}(a,b))=1-\frac{TLB(\mathbf{d}(a,b))}{TUB(\mathbf{d}(a,b))}=1-2\frac{1-\frac{1}{a+b}}{1-\frac{3}{2a+b}+\sqrt{1+\frac{ab}{2a+b}+\frac{2b-4a+1}{2a+b}}}.
\end{equation}

Introduce the new variable, $n:=a+b$, and exclude $b$ from
(\ref{eq_TErr_ab}):
\begin{equation}\label{eq_TErr_ab2}
TErr(\mathbf{d}(a,n-a))=1-2\frac{1-\frac{1}{n}}{1-\frac{3}{n+a}+\sqrt{1+\frac{a(n-a)}{(n+a)^2}+\frac{2n+1-6a}{(n+a)^2}}}.
\end{equation}

Find the maximum by $a$ of the right-hand side of
(\ref{eq_TErr_ab2}) assuming $n$ constant. From the first-order
conditions, the extremal $a$ satisfies the equation
$$6\sqrt{3an+n^2-6a+2n+1}+n^2-3an+6a=10n+2.$$
Omitting routine calculations, we conclude that for $n\ge2$ the
maximum of (\ref{eq_TErr_ab2}) is attained at
$a(n):=\frac{n-8}{3}+2\sqrt{2}$, and
$$TErr(\mathbf{d}(a(n),n-a(n)))=
1-\frac{8(n-2+3\sqrt{2})(1-\frac{1}{n})}{4n-17+6\sqrt{2}+3\sqrt{2n^2-8n+6n\sqrt{2}+17-12\sqrt{2}}}.$$

Standard analysis shows that $TErr(\mathbf{d}(a(n),n-a(n)))$ is a
monotone function of $n$. Therefore, for any $a,b \in \mathbb{N}$
$$TErr(\mathbf{d}(a,b))\le
\lim_{n\rightarrow+\infty}TErr(\mathbf{d}(a(n),n-a(n)))=\frac{3\sqrt{2}-4}{3\sqrt{2}+4}.$$
\end{proof}
\end{theorem}

\section{Conclusion}
Above we suggested a lower bound for the terminal distance spectral
radius of a tree with the given degree sequence and showed it to be
within $3\%$ from the terminal distance spectral radius of the
$BFS$-tree. This means that Conjecture \ref{conj2}, which says that
the $BFS$-tree has minimal $TDSR$ among all tree with the given
degree sequence, is valid at least up to $3\%$ (and typically even
more precisely, as shown in Fig.~\ref{figerror}).

However, our proof of this $3\%$ error relies on the Conjecture
\ref{conj3}, which guesses the shape of trees, which give the
maximum error. Although being pretty natural and well-grounded
numerically, this conjecture is still an open issue in our analysis.

It is of interest to discuss how much the presented technique can be
applied to estimate the precision of Conjecture \ref{conj1}, which
states a similar minimum property of the $BFS$-tree with respect to
the distance spectral radius ($DSR$).

In Fig.~\ref{figerrorWI} the relative error
$$Err(\mathbf{d}):=\frac{UB(\mathbf{d}(a,b))-LB(\mathbf{d}(a,b))}{UB(\mathbf{d}(a,b))},$$
where $LB$ and $UB$ are defined in (\ref{eq_main_ineq_dsr}), is
presented for all trees of order $N\le 23$. We see that the error of
$LB$ is, in average, twice as big as that of $TLB$ (see
Fig.~\ref{figerror}). From Fig.~\ref{figerrorWI} one might
conjecture that $Err$ does not exceed $6\%$, but careful analysis,
similar to that performed in Lemma \ref{lemma_TDSR_ab}, needs to be
performed to justify this conjecture.

An observation, which may help, is that degree sequences giving the
maximum $Err$ (and, hence, forming the upper envelope of the point
set in Fig.~\ref{figerrorWI}) have the form
$(\underbrace{1,...,1}_{d\textrm{ times}}, 2,...,2,d).$ They
correspond to, the so-called, \emph{starlike trees}, whose distance
spectral properties are studied in detail by Stevanovi\'{c} and
Ili\'{c} \cite{StevanovicIlic2010}.

\vspace{4mm}

\begin{center}
\includegraphics[width=7cm,keepaspectratio]{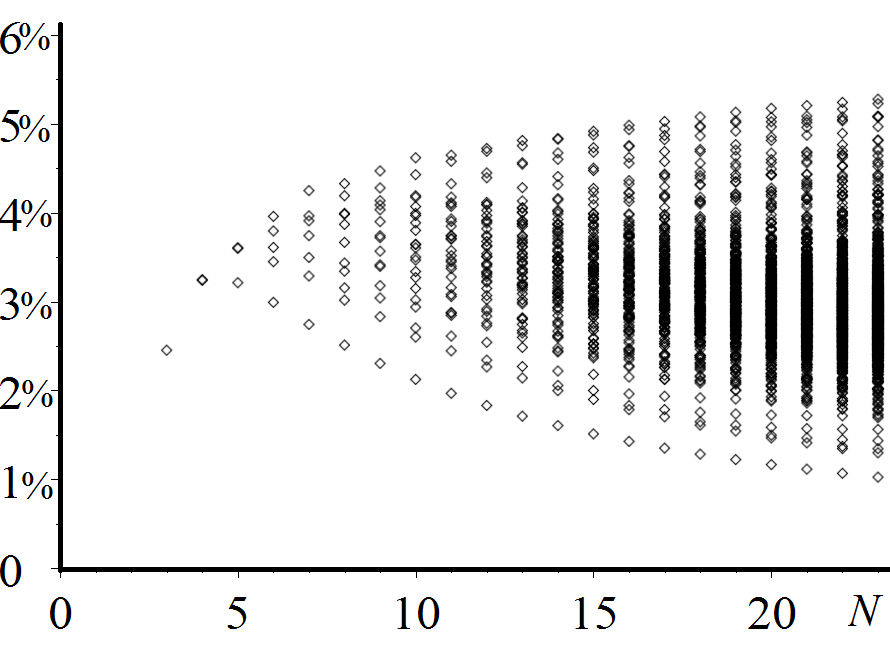}
\label{figerrorWI}
\end{center}
\baselineskip=0.20in

\noindent {\bf Fig. 4.} Relative error $Err$ of the lower bound for
$DSR$ vs the tree order $N$.

\baselineskip=0.30in

\vspace{3mm}

\vspace{8mm}

\baselineskip=0.25in

\noindent {\bf Acknowledgements.}

\vspace{1mm}

This research is supported by the grant of the Russian Foundation
for Basic Research, No 13-07-00389.

\vspace{8mm}

\baselineskip=0.25in

\section{Appendix. Proof of Theorem \ref{theorem_bfs_optimal}}
First we recall the notions of a directed tree and subordinate group
weight from \cite{Goubko2015}.

\begin{definition}
A vertex-weighted \emph{directed tree} is a connected directed graph
where each vertex except the \emph{root} has the sole outbound arc,
the root has no outbound arcs, and each vertex has non-negative
weight assigned to it.
\end{definition}

An arbitrary vertex-weighted tree $T$ consisting of more than two
vertices can be transformed into a directed tree $T_r$ by choosing
an \textbf{internal vertex} $r\in M(T)$ as a root, and replacing all
its edges with arcs directed towards the root. Let
$\mathcal{WR}(\mu, d)$ stand for all directed trees obtained from
$\mathcal{WT}(\mu,d)$ in this way.

If at Step $i$ of the Huffman algorithm (see Section
\ref{section_preliminaries}) directed arcs towards the vertex $m_i$
are added instead of undirected edges, the algorithm builds a
\emph{directed Huffman tree}. Let $\mathcal{RH}(\mu, d)$ be the set
of directed Huffman trees for a generating tuple $\langle\mu,
d\rangle$.

\begin{definition}
For an arbitrary vertex $v\in V(T)$ of the directed tree $T$ define
its \emph{subordinate group} $g_T(v)\subseteq V(T)$ as the set of
vertices having the directed path to the vertex $v$ in the tree $T$
\emph{(}the vertex $v$ itself belongs to $g_T(v)$\emph{)}. The
\emph{weight $f_T(v)$ of the subordinate group} $g_T(v)$ is defined
as the total vertex weight of the group: $f_T(v):=\sum_{u\in g_T(v)}
\mu_T(u)$.
\end{definition}

\begin{definition}\label{def_vector}\sloppy
For a directed tree $T\in \mathcal{WR}(\mu,d)$ with root $r$ define
a \emph{vector} $\mathbf{f}(T):=(f_T(m))_{m\in
M(T)\backslash\{r\}}^\uparrow$ \emph{of subordinate groups'
weights}.
\end{definition}

If some tree $T\in \mathcal{WT}(\mu,d)$ is transformed into a
directed tree $T_r\in \mathcal{WR}(\mu,d)$ by choosing a root $r$,
its Wiener index can be written as \cite{Klavzar1997}:
\begin{equation}\label{eq_VWWI_directed}
VWWI(T)=VWWI(T_r)=\sum_{v\in V(T)\backslash
\{r\}}f_{T_r}(v)(\bar{\mu}-f_{T_r}(v))=\sum_{v\in V(T)\backslash
\{r\}}\chi(f_{T_r}(v)),
\end{equation}
where $\chi(x):=x(\bar{\mu}-x)$.

\begin{lemma}\label{lemma_Huffman_monotonicity}\textbf{\emph{\cite{Goubko2015}}}
If weights are degree-monotone in $\langle\mu,d\rangle$, then for
any $H\in \mathcal{RH}(\mu,d)$
\begin{equation}\label{eq_weights_monotone}
[vm,v'm'\in E(H), m \neq m', f_H(v)<f_H(v')] \Rightarrow f_H(m) <
f_H(m').
\end{equation}
\end{lemma}

\begin{definition}\textbf{\emph{\cite{Marshall79,Zhang2008W}}}
A non-negative sequence $\mathbf{x} = (x_1, ..., x_p)$, $p \in
\mathbb{N}$, \emph{weakly majorizes} a non-negative sequence
$\mathbf{y} = (y_1, ..., y_p)$ \emph{(}denote it
$\mathbf{y}\preceq^w\mathbf{x}$ or
$\mathbf{x}\succeq^w\mathbf{y}$\emph{)} if
$$\sum_{i=1}^k x_{[i]} \le \sum_{i=1}^k y_{[i]} \text{ for all }k=1,...,p.$$
If $\mathbf{x}^\uparrow \neq \mathbf{y}^\uparrow$, then $\mathbf{x}$
is said to \emph{strictly weakly majorize} $\mathbf{y}$ $($denote it
$\mathbf{y}\prec^w \mathbf{x}$ or $\mathbf{x}\succ^w \mathbf{y}$$)$.
\end{definition}

We will need the following properties of weak majorization.

\begin{lemma}\label{lemma_Zhang_b} \textbf{\emph{\cite{Marshall79,Zhang2008W}}}
Consider a positive number $b>0$ and two non-negative sequences,
$\mathbf{x} = (x_1, ..., x_k, y_1, ..., y_l)$ and $\mathbf{y} = (x_1
+ b, ..., x_k + b, y_1 - b, ..., y_l - b)$, such that $0 \le k \le
l$. If $x_i \ge y_i$ for $i = 1, ..., k$, then $\mathbf{x}\prec^w
\mathbf{y}$.
\end{lemma}

\begin{lemma}\label{lemma_Zhang_xy} \textbf{\emph{\cite{Marshall79,Zhang2008W}}}
If $\mathbf{x}\preceq^w \mathbf{y}$ and $\mathbf{x}' \preceq^w
\mathbf{y'}$, then $(\mathbf{x},\mathbf{x'})\preceq^w
(\mathbf{y},\mathbf{y'})$, where $(\mathbf{x},\mathbf{x'})$ is
concatenation of sequences $\mathbf{x}$ and $\mathbf{x'}$. Moreover,
if $\mathbf{x}' \prec^w \mathbf{y'}$, then
$(\mathbf{x},\mathbf{x'})\prec^w (\mathbf{y},\mathbf{y'})$.
\end{lemma}

\begin{lemma}\label{lemma_Zhang_concave} \textbf{\emph{\cite{Marshall79,Zhang2008W}}}\sloppy
If $\chi(x)$ is an increasing concave function and
$(x_1,...,x_p)\preceq^w (y_1,...,y_p)$, then
$\sum_{i=1}^p\chi(x_i)\ge \sum_{i=1}^p\chi(y_i)$, and equality is
possible only when $(x_1,...,x_p)^\uparrow=(y_1,...,y_p)^\uparrow$.
\end{lemma}

\begin{theorem}\label{theorem_major}\textbf{\emph{\cite{Goubko2015}}}
If weights are degree-monotone in a generating tuple $\langle\mu,
d\rangle$ and $H\in \mathcal{RH}(\mu, d)$, then for any directed
tree $T \in \mathcal{WR}(\mu, d)$ $\mathbf{f}(H)\succeq^w
\mathbf{f}(T)$.
\end{theorem}

\begin{lemma}\label{lemma_equal_weights}
Consider a degree sequence $\mathbf{d}$ with $n$ elements being
equal to unity and a positive sequence $\mu=(\mu_1,...,\mu_n)$ such
that $\mu_i>\mu_j$ for some $i$ and $j$. Define a generating tuple
$\langle\mu',d\rangle$, which differs from $\langle\mu,d\rangle$
only with $i$-th and $j$-th weight components, namely,
$\mu'_i-\varepsilon=\mu'_j+\varepsilon=\frac{1}{2}(\mu_i+\mu_j)$,
where $0<\varepsilon<\delta:=\frac{1}{2}(\mu_i-\mu_j)$. If $H\in
\mathcal{RH}(\mu,d)$ and $H'\in \mathcal{RH}(\mu',d)$ are some
directed Huffman trees, then $\mathbf{f}(H)\succeq^w\mathbf{f}(H')$.
\begin{proof}
Suppose vertices $u, v\in W(\mu,d)$ have weights $\mu_i$ and $\mu_j$
respectively. Let $(u, u_1, ..., u_k, m)$ and $(v, v_1, ..., v_l,
m)$ be the disjoint paths in the directed Huffman tree $H'\in
\mathcal{RH}(\mu', d)$ from vertices $u$ and $v$ to some vertex $m
\in M(H')$, where $k,l\ge 0$. Since
$f_{H'}(u)=\mu'_i,f_{H'}(v)=\mu'_j$, and $\mu'_i>\mu'_j$, it follows
immediately from Lemma \ref{lemma_Huffman_monotonicity} that $k\le
l$ and $f_{H'}(u_i)>f_{H'}(v_i)$ for $i=1,...k$.

Consider a directed weighted tree $T$ obtained from $H'$ by changing
the weights of vertices $u,v$ to $\mu_i$ and $\mu_j$ respectively.
In the tree $T$ weights of groups subordinated to the vertices $u_1,
..., u_k$ increase by $b:=\delta-\varepsilon>0$ (i.e., $f_{T}(u_i) =
f_{H'}(u_i) + b, i = 1, ..., k$), weights of the groups subordinated
to the vertices $v_1, ..., v_l$ decrease by $b$ (i.e., $f_{T}(v_i) =
f_{H'}(v_i) - b, i = 1, ..., l$), weights of all other vertices
(including $m$) do not change.

Therefore, by Lemma \ref{lemma_Zhang_b},
$$\mathbf{y} := (f_{T}(u_1), ..., f_{T}(u_k), f_{T}(v_1), ..., f_{T}(v_l))= $$
$$=(f_{H'}(u_1) + b, ..., f_{H'}(u_k) + b, f_{H'}(v_1) - b, ..., f_{H'}(v_l) - b) \succ^w $$
$$\succ^w (f_{H'}(u_1), ..., f_{H'}(u_k), f_{H'}(v_1), ..., f_{H'}(v_l)) =: \mathbf{x}.$$

If one denotes with $\mathbf{z}$ the sequence of (unchanged) weights
of groups subordinated to all other internal vertices of $T$
distinct from the root, then, by Lemma \ref{lemma_Zhang_xy},
$\mathbf{f}(T)=(\mathbf{y},\mathbf{z})\succ^w
(\mathbf{x},\mathbf{z})=\mathbf{f}(H')$.

It is clear that $T\in \mathcal{WR}(\mu,d)$. Since $H\in
\mathcal{RH}(\mu,d)$, from Theorem \ref{theorem_major} we know that
$\mathbf{f}(H)\succeq^w \mathbf{f}(T)$ and, consequently,
$\mathbf{f}(H)\succeq^w \mathbf{f}(H')$.
\end{proof}
\end{lemma}

\begin{corollary}\label{cor_simplex}
Under the conditions of Lemma \ref{lemma_equal_weights}
$TVWWI(H)<TVWWI(H')$.
\begin{proof}
Let $r$ and $r'$ be roots of trees $H$ and $H'$ respectively.
According to (\ref{eq_VWWI_directed}),
$$TVWWI(H)=\sum_{k=1}^n\chi(\mu_k)+\sum_{v\in
M(\mu, d)\backslash \{r\}}\chi(f_{H}(v)),\textrm{ where }
\chi(x)=x(\bar{\mu}-x),$$ and
$$TVWWI(H')=\sum_{k=1}^n\chi(\mu_k')+\sum_{v\in
M(\mu, d)\backslash \{r'\}}\chi(f_{H'}(v)).$$

From Lemma \ref{lemma_equal_weights} we know that
$\mathbf{f}(H)\succeq^w \mathbf{f}(H')$. Therefore, by Lemma
\ref{lemma_Zhang_concave},
$$\sum_{v\in M(\mu, d)\backslash
\{r\}}\chi(f_{H}(v))< \sum_{v\in M(\mu', d)\backslash
\{r'\}}\chi(f_{H'}(v)).$$

Since $\chi(x)$ is concave,
$$\chi(\mu_i)+\chi(\mu_j)<\chi(\mu_i-\delta+\varepsilon)+\chi(\mu_j+\delta-\varepsilon).$$
Other vertex weights do not change, so $TVWWI(H)<TVWWI(H')$.
\end{proof}
\end{corollary}

\begin{lemma}\label{lemma_simplex_sphere}
Consider a degree sequence $\mathbf{d}$ with $n$ elements being
equal to unity and a vector $\mu\in S_n^+$. If $\mu_i>\mu_j$ for
some $i$ and $j$, then such a vector $\nu\in S_n^+$ exists that
$TVWWI(H'')<TVWWI(H)$, where $H''\in \mathcal{H}(\mu,d)$ and $H\in
\mathcal{H}(\nu,d)$ are Huffman trees.
\begin{proof}
Recall a weight vector $\mu'$ from Lemma \ref{lemma_equal_weights}
and define the vector $\nu\in S_n^+$ as
$$\nu_i=\sqrt{\frac{\mu_i^2+\mu_j^2}{2}+(\mu_i+\mu_j)\varepsilon},$$
$$\nu_i=\sqrt{\frac{\mu_i^2+\mu_j^2}{2}-(\mu_i+\mu_j)\varepsilon},$$
$$\nu_k=\mu_k, k\neq i,j.$$
Interrelation between $\mu$, $\mu'$, and $\nu$ is shown in
Fig.~\ref{figure_mu_nu}.

Let $\mu_H(u)=\mu_i$, $\mu_H(v)=\mu_j$ for some $u,v\in V(H)$ and
consider a vertex-weighted tree $T$ obtained from $H$ by changing
weights of $u$ and $v$ to $\mu_i'$ and $\mu_j'$ respectively. From
Fig.~2 is is clear that $\nu_i>\mu_i'$, $\nu_j>\mu_j'$. Since
$RD(\cdot)$ is off-diagonal positive, $TVWWI(\cdot)$ is strictly
monotone in weights, so $TVWWI(T)<TVWWI(H)$. Consider a tree $H'\in
\mathcal{H}(\mu',d)$. By Note \ref{note_tvwwi_monotone} we have
$TVWWI(H')\le TVWWI(T)$. From Corollary \ref{cor_simplex} we also
know that $TVWWI(H'')<TVWWI(H')$. Summarizing these inequalities
obtain the statement of the lemma.
\end{proof}
\end{lemma}

\vspace{4mm}

\begin{center}
\includegraphics[height=8cm,keepaspectratio]{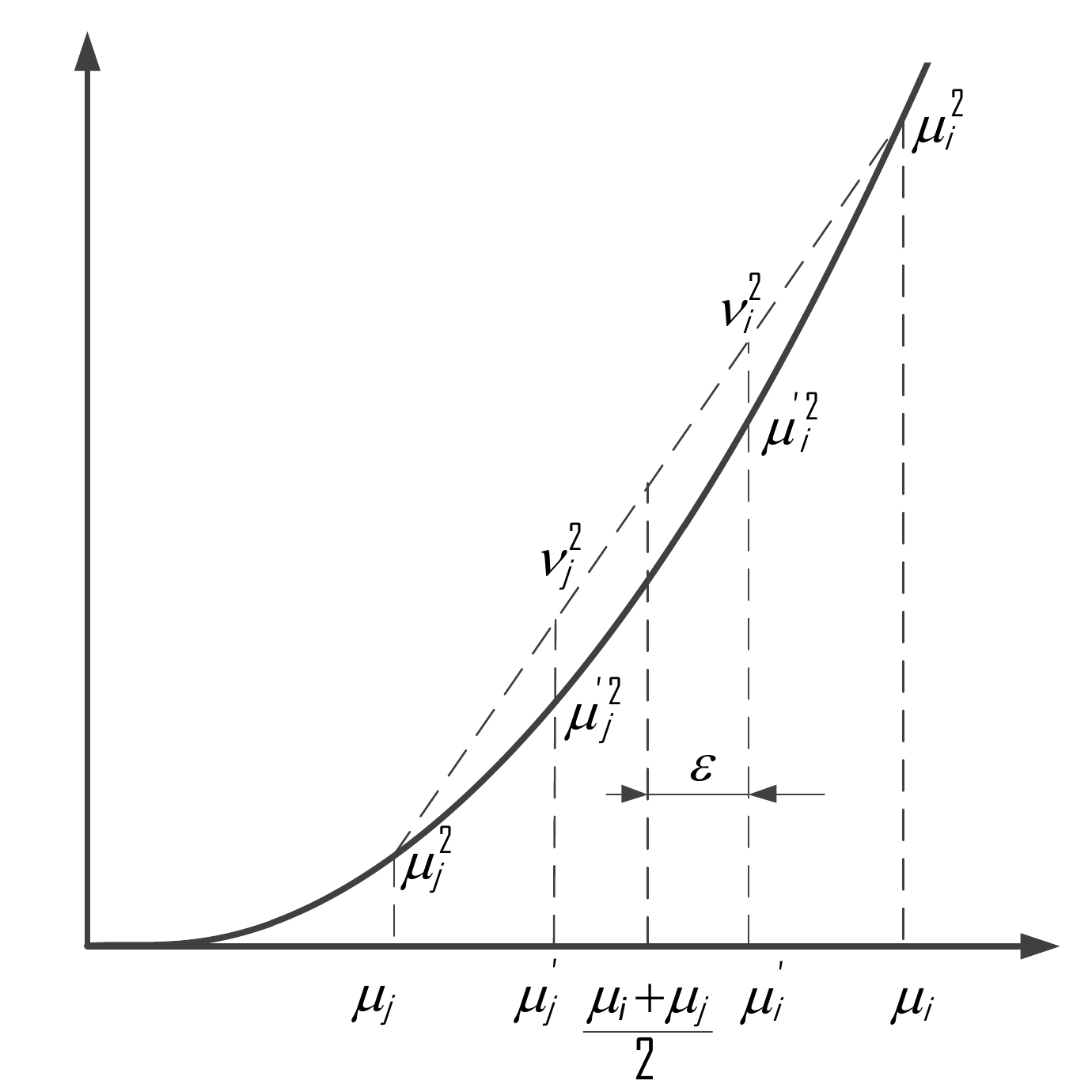}
\label{figure_mu_nu}
\end{center}
\baselineskip=0.20in

\noindent {\bf Fig. 5.} Interrelation between $\mu$, $\mu'$,and
$\nu$.

\baselineskip=0.30in

\vspace{3mm}

Lemma \ref{lemma_simplex_sphere} says that when positive weights for
pendent vertices are picked from the unit sphere, $TVWWI(\cdot)$ of
a tree with the degree sequence achieves its maximum when all vertex
weights are equal to each other, i.e., $\mu_i=1/\sqrt{n}$ for all
$i=1,...,n$.

The vector
$\mathbf{e}:=(\frac{1}{\sqrt{n}},...,\frac{1}{\sqrt{n}})^T$ of equal
weights is the only positive vector on the unit sphere, to which
Lemma \ref{lemma_simplex_sphere} cannot be applied, and, since for
an arbitrary graph $G$ the identity holds $\mathbf{e}^T
RD(G)\mathbf{e}=\mathbf{1}^T RD(G)\mathbf{1}/n=2TWI(G)/n$, the the
statement of Theorem~\ref{theorem_bfs_optimal} follows immediately.

\end{document}